\documentclass[11pt]{article}
\usepackage[margin=1in]{geometry} 
\usepackage{amssymb,amsfonts,amsmath,bbm,mathrsfs,stmaryrd}
\usepackage{xcolor}
\usepackage{url}

\usepackage[colorlinks,
             linkcolor=black!75!red,
             citecolor=blue,
             pdftitle={Cosemisimple Hopf algebras are faithfully flat over Hopf subalgebras},
             pdfauthor={Alexandru Chirvasitu},
             pdfsubject={rings and algebras},
             pdfkeywords={ra.rings and algebras, ct.category theory, rt.representation theory, oa.operator algebras},
             pdfproducer={pdfLaTeX},
             pdfpagemode=None,
             bookmarksopen=true
             bookmarksnumbered=true]{hyperref}

\usepackage{tikz}
\usetikzlibrary{arrows,decorations.pathreplacing,decorations.markings,shapes.geometric,through,fit,shapes.symbols,positioning}

\usepackage[amsmath,thmmarks,hyperref]{ntheorem}
\usepackage{cleveref}

\crefname{section}{Section}{Sections}
\crefformat{section}{#2Section~#1#3} 
\Crefformat{section}{#2Section~#1#3} 

\crefname{subsection}{\S}{\S\S}
\crefformat{subsection}{#2\S#1#3} 
\Crefformat{subsection}{#2\S#1#3} 

\theoremstyle{change}

\newtheorem{lemma}{Lemma}[section]
\newtheorem{proposition}[lemma]{Proposition}
\newtheorem{corollary}[lemma]{Corollary}
\newtheorem{theorem}[lemma]{Theorem}

\theoremstyle{nonumberplain}
\newtheorem{theoremN}{Theorem}

\theoremstyle{change}
\theorembodyfont{\upshape}
\theoremsymbol{\ensuremath{\blacklozenge}}

\newtheorem{definition}[lemma]{Definition}

\newtheorem{remark}[lemma]{Remark}

\crefname{definition}{definition}{definitions}
\crefformat{definition}{#2definition~#1#3} 
\Crefformat{definition}{#2Definition~#1#3} 

\crefname{lemma}{lemma}{lemmas}
\crefformat{lemma}{#2lemma~#1#3} 
\Crefformat{lemma}{#2Lemma~#1#3} 

\crefname{proposition}{proposition}{propositions}
\crefformat{proposition}{#2proposition~#1#3} 
\Crefformat{proposition}{#2Proposition~#1#3} 

\crefname{example}{example}{examples}
\crefformat{example}{#2example~#1#3} 
\Crefformat{example}{#2Example~#1#3} 

\crefname{remark}{remark}{remarks}
\crefformat{remark}{#2remark~#1#3} 
\Crefformat{remark}{#2Remark~#1#3} 

\crefname{corollary}{corollary}{corollaries}
\crefformat{corollary}{#2corollary~#1#3} 
\Crefformat{corollary}{#2Corollary~#1#3} 

\crefname{theorem}{theorem}{theorems}
\crefformat{theorem}{#2theorem~#1#3} 
\Crefformat{theorem}{#2Theorem~#1#3} 

\crefname{equation}{}{}
\crefformat{equation}{(#2#1#3)} 
\Crefformat{equation}{(#2#1#3)}

\theoremstyle{nonumberplain}
\theoremsymbol{\ensuremath{\blacksquare}}

\newtheorem{proof}{Proof}
\newtheorem{sketch}{Sketch of proof}
\newtheorem{proof of Aproj}{Proof of \Cref{prop.Aproj means proj}}
\newtheorem{proof of thCQG}{Proof of \Cref{th.CQG}}

\DeclareMathOperator{\id}{id}
\DeclareMathOperator{\tr}{tr}

\newcommand\Ff{{\mathcal F}}

\newcommand\Mm{{\mathcal M}}

\newcommand\Oo{{\mathcal O}}

\newcommand\CC{{\mathbb C}}

\newcommand{\ol}{\overline}

\newcommand\tofrom{\longleftrightarrow}

\newcommand\longto{\longrightarrow}

\newcommand\verylongmapsto[1]{\!{\tikz[anchor=base,auto,baseline=(a.base)] \draw[|-angle 90] node (a) {\phantom{A}} ++(1cm,0) node (b) {\phantom{A}\!\!\!} (a.mid) -- node{$\scriptstyle #1$} (b.mid);}}

\DeclareMathOperator{\ad}{{\rm ad}}

\DeclareMathOperator{\Hom}{Hom}

\DeclareMathOperator{\spec}{spec}

\DeclareMathOperator{\Functors}{Functors}

\DeclareMathOperator*{\revarinjlim}{{\tikz[baseline=(lim.base)] \draw[->] node[inner sep = 0pt] (lim) {$\lim$} (lim.south west) ++(1pt,-3pt) coordinate (a) (lim.south east) ++(-1pt,-3pt) coordinate (b) (a) -- (b) ;}}
\renewcommand\varinjlim\revarinjlim

\DeclareMathOperator*{\revarprojlim}{{\tikz[baseline=(lim.base)] \draw[<-] node[inner sep = 0pt] (lim) {$\lim$} (lim.south west) ++(1pt,-3pt) coordinate (a) (lim.south east) ++(-1pt,-3pt) coordinate (b) (a) -- (b) ;}}
\renewcommand\varprojlim\revarprojlim

\newcommand\Fun\Functors

\newcommand{\define}[1]{{\em #1}}

\newcommand{\cat}[1]{\textsc{#1}}

\newcommand{\qedhere}{\mbox{}\hfill\ensuremath{\blacksquare}}

\renewcommand{\square}{\mathrel{\Box}}

\title{Cosemisimple Hopf algebras are faithfully flat over Hopf subalgebras}
\author{Alexandru Chirvasitu\footnote{UC Berkeley, \url{chirvasitua@math.berkeley.edu}}}

\begin{document}
\maketitle

\begin{abstract}
	The question of whether or not a Hopf algebra $H$ is faithfully flat over a Hopf subalgebra $A$ has received positive answers in several particular cases: when $H$ (or more generally, just $A$) is commutative, or cocommutative, or pointed, or when $K$ contains the coradical of $H$. We prove the statement in the title, adding the class of cosemisimple Hopf algebras to those known to be faithfully flat over all Hopf subalgebras. We also show that the third term of the resulting ``exact sequence'' $A\to H\to C$ is always a cosemisimple coalgebra, and that the expectation $H\to A$ is positive when $H$ is a CQG algebra. 	
\end{abstract}

\noindent {\em Keywords: cosemisimple Hopf algebra, CQG algebra, faithfully flat, right coideal subalgebra, quotient left module coalgebra, expectation}

\tableofcontents


\section*{Introduction}

	The issue of faithful flatness of a Hopf algebra (always over a field) over its Hopf subalgebras arises quite naturally in several ways. One direction is via the so-called Kaplansky conjecture (\cite{MR0435126}), which initially asked whether or not Hopf algebras are free over Hopf subalgebras (as an analogue to the Lagrange theorem for finite groups). The answer was known to be negative, with a counterexample having appeared in \cite{MR0360610}, but it is true in certain particular cases: using the notations in the abstract, $H$ is free over $A$ whenever $H$ is finite dimensional (Nichols-Zoeller Theorem, \cite[Theorem 3.1.5]{MR1243637}), or pointed (\cite{MR0437582}), or $A$ contains the coradical of $H$ (\cite[Corollary 2.3]{MR0437581}). 
	
	Montgomery then naturally asks whether one can get a positive result by requiring only faithful flatness of a Hopf algebra over an arbitrary Hopf subalgebra (\cite[Question 3.5.4]{MR1243637}). Again, this turns out not to work in general (see \cite{MR1761130} and also \cite{MR2584972}, where the same problem is considered in the context of whether or not epimorphisms of Hopf algebras are surjective), but one has positive results in several important cases, such as that when $A$ is commutative (\cite[Proposition 3.12]{MR1954457}), or $H$ is cocommutative (\cite[Theorem 3.2]{MR0321963}, which also takes care of the case when $H$ is commutative). The most recent version of the question, asked in \cite{MR1761130}, seems to be whether or not a Hopf algebra with bijective antipode is faithfully flat over Hopf subalgebras with bijective antipode.   
	
	Another way to get to the faithful flatness issue is via the problem of constructiing quotients of affine group schemes. We recall briefly how this goes.  
	 
	Let $A\to H$ be an inclusion of commutative Hopf algebras; in scheme language, $A$ and $H$ are affine groups, and the inclusion means that $\spec(A)$ is a quotient group scheme of $\spec{H}$. The Hopf algebraic analogue of the kernel of this epimorphism is the quotient Hopf algebra $\pi:H\to C=H/HA^+$, where $A^+$ stands for the kernel of the counit of $A$. The map $\pi$ is then \define{normal}, in the sense of \cite[Definition 1.1.5]{Andruskiewitsch1995}:
	\[
		\cat{LKer}(\pi)=\{a\in A\ |\ (\pi\otimes\id)\circ\Delta(a)=1_C\otimes a\}
	\]  
equals its counterpart
	\[
		\cat{RKer}(\pi)=\{a\in A\ |\ (\id\otimes\pi)\circ\Delta(a)=a\otimes 1_C\}. 
	\]
This means precisely that $\spec(C)$ is a normal affine subgroup scheme of $\spec(A)$ (\cite[Lemma 5.1]{MR0321963}). This gives a map $A\mapsto C$ from quotient affine group schemes of $H$ to normal subgroup schemes. One naturally suspects that this is probably a bijective correspondence, and this is indeed true (see \cite[Theorem 4.3]{MR0321963} and also \cite[III $\S$3 7.2]{MR0302656}). In Takeuchi's paper faithful flatness is crucial in proving half of this result, namely the injectivity of the map $A\mapsto C$: one recovers $A$ as $\cat{LKer}(\pi)$.

	Many of the technical arguments and constructions appearing in this context go through in the non-commutative setting, so one might naturally be led to the faithful flatness issue by trying to mimic the algebraic group theory in a more general setting, where Hopf algebras are viewed as function algebras on a ``quantum'' group. This is, for example, the point of view taken in the by now very rich and fruitful theory of compact quantum groups, first introduced and studied by Woronowicz: the main characters are certain $C^*$ algebras $A$ with a comultiplication $A\to A\underline{\otimes} A$ (minimal $C^*$ tensor product), imitating the algebras of continuous functions on a compact group (we refer the reader to \cite[Chapter 11]{MR1492989} or Woronowicz's landmark papers \cite{MR901157,MR943923} for details). 
	
	These objects are not quite Hopf algebras, but for any compact quantum group $A$ as above, one can introduce a genuine Hopf algebra $\mathcal A$, imitating the algebra of \define{representative} functions on a compact group (i.e. linear span of matrix coefficients of finite dimensional unitary representations), and which contains all the relevant information on the representation theory of the quantum group in question. The abstract properties of such Hopf ($*-$)algebras have been axiomatized, and they are usually referred to as CQG algebras (see \cite[11.3]{MR1492989} or the original paper \cite{MR1310296}, where the term was coined). They are always cosemisimple (as an analogue of Peter-Weyl theory for representations of compact groups), which is why we hope that despite the seemingly restrictive hypothesis of cosemisimplicity, the results in the present paper might be useful apart from any intrinsic interest, at least in dealing with Hopf algebraic issues arising in the context of compact quantum groups. 
			
	We now describe the contents of the paper. 
	
	In the first section we introduce the conventions and notations to be used throughout the rest of the paper, and also develop the tools needed to prove the main results. In \Cref{subse.descent} we set up a Galois correspondence between the sets of right coideal subalgebras of a Hopf algebra $H$ and the set of quotient left module coalgebras of $H$. We then recall basic results on categories of objects imitating Sweedler's Hopf modules: These have both a module and a comodule structure, one of them over a Hopf algebra $H$, and the other one over a right coideal subalgebra or a quotient left module coalgebra of $H$. These categories are used extensively in the subsequent discussion. 
	
	\Cref{se.main} is devoted to the main results. We provide sufficient conditions for faithful flatness over Hopf subalgebras in \Cref{th.main,cor.main}. We also investigate the case of cosemisimple $H$ further, proving in \Cref{th.mainbis} that for any Hopf subalgebra $A$, the quotient left $H$-module coalgebra $C=H/HA^+$ is always cosemisimple. This quotient is the third term of the ``exact sequence'' which completes the inclusion $A\to H$, and the question of whether or not $C$ is cosemisimple arises naturally in the course of the proof of \Cref{th.main}, which shows immediately that it is true when $HA^+$ happens to be an ideal (both left \define{and} right).   
	
	Finally, \Cref{se.CQG} we show that when the ambient Hopf algebra $H$ is CQG, the ``expectation'' $H\to A$ that plays a crucial role in the preceding section is positive. In the course of the proof we use a sort of ``$A$-relative'' Fourier transform from $H$ to $C^*$ (whereas ordinary Fourier transforms, as in, say, \cite{MR1059324}, are roughly speaking more like maps from $H$ to the dual $H^*$). This construction has some of the familiar properties from harmonic analysis, such as intertwining products and ``convolution products'' (\Cref{prop.F} 1.), playing well with $*$ structures (\Cref{prop.F} 2.), and satisfying a Plancherel-type condition (\Cref{rem.Planch}).

\subsection*{Acknowledgements}

I would like to thank Akira Masuoka and Issan Patri for useful references, conversations and comments on the contents of the paper. I am also grateful to the anonymous referee for bringing to my attention references \cite{MR688207,MR1098982}, and for numerous other suggestions that have helped considerably with improving the manuscript.


\section{Preliminaries}\label{se.prelim}

	In this section we make the preparations necessary to prove the main results. Throughout, we work over a fixed field $k$, so (co, bi, Hopf)algebra means (co, bi, Hopf)algebra over $k$, etc. The reader should feel free to assume $k$ to be algebraically closed whenever convenient, as most results are invariant under scalar extension. In \Cref{se.CQG} we specialize to characteristic zero. 
	
	We assume basic familiarity with coalgebra and Hopf algebra theory, for example as presented in \cite{MR1243637}. We will make brief use of the notion of coring over a (not necessarily commutative) $k$-algebra; we refer to \cite{MR2012570} for basic properties and results. 
	
	The notations are standard: $\Delta_C$ and $\varepsilon_C$ stand for comultiplication and counit of the coalgebra $C$ respectively, and we will allow ourselves to drop the subscript when it is clear which coalgebra is being discussed. Similarly, $S_H$ or $S$ stands for the antipode of the Hopf algebra $H$, $1_A$ (or just $1$) will be the unit of the algebra $A$, etc. Sweedler notation for comultiplication is used throughout, as in $\Delta(h)=h_1\otimes h_2$, as well as for left or right coactions: if $\rho:N\to N\otimes C$ ($\rho:N\to C\otimes N$) is a right (left) $C$-comodule structure, we write $n_0\otimes n_1$ ($n_{\langle -1\rangle}\otimes n_{\langle 0\rangle}$) for $\rho(n)$. We sometimes adorn the indices with parentheses, as in $\Delta(c)=c_{(1)}\otimes c_{(2)}$.    
	
	We will also be working extensively with categories of (co)modules over (co)algebras, as well as categories of objects admitting both a module and a comodule structure, compatible in some sense that will be made precise below (see \Cref{subse.descent}). These categories are always denoted by the letter $\Mm$, with left (right) module structures appearing as left (right) subscripts, and left (right) comodule structures appearing as left (right) superscripts. All such categories are abelian (and in fact Grothendieck), and the forgetful functor from each of them to vector spaces is exact. The one exception from this notational convention is the category of $k$-vector spaces, which we simply call $\cat{Vec}$.
	
	Recall that the category $\Mm^H_f$ of finite dimensional right comodules over a Hopf algebra is monoidal left rigid: every object $V$ has a left dual $V^*$ (at the level of vector spaces it is just the usual dual vector space), and one has adjunctions $(\otimes V,\otimes V^*)$ and $(V^*\otimes,V\otimes)$ (the left hand member of the pair is the left adjoint) on $\Mm^H_f$.	
		
	We also use the correspondence between subcoalgebras of a Hopf algebra $H$ and finite dimensional (right) comodules over $H$: for such a comodule $V$, there is a smallest subcoalgebra $D=\cat{coalg}(V)\le H$ such that the structure map $V\to V\otimes H$ factors through $V\to V\otimes D$. Conversely, if $D\le H$ is a simple subcoalgebra, then we denote by $V_D$ the simple right $D$-comodule, viewed as a right $H$-comodule. Then, for simple subcoalgebras $D,E\le H$, the product $ED$ will be precisely $\cat{coalg}(V_E\otimes V_D)$, while $S(D)$ is $\cat{coalg}(V^*)$. 
	
	For a coalgebra $C$, the symbol $\widehat{C}$ denotes the set of isomorphism classes of simple (right, unless specified otherwise) $C$-comodules.


\subsection{Descent data and adjunctions}\label{subse.descent}

We will be dealing with the kind of situation studied extensively in \cite{MR549940}: $H$ will be a Hopf algebra, and for most of this section (and in fact the paper), $\iota:A\to H$ will be a right coideal subalgebra, while $\pi: H\to C$ will be a quotient left $H$-module coalgebra. Recall that this means that $A$ is a right coideal of $H$ (i.e. $\Delta_H(A)\le A\otimes H$) as well as a subalgebra, and so the induced map $A\to A\otimes H$ is an algebra map; similarly, $C$ is the quotient of $H$ by a left ideal as well as a coalgebra, and the induced map $H\otimes C\to C$ is supposed to be a coalgebra map.   
	
Given a coalgebra map $\pi:H\to C$, we write $\ol h$ for $\pi(h)$, $h\in H$. In this situation, $H$ will naturally be both a left and a right $C$-comodule (via the structure maps $(\pi\otimes\id)\circ\Delta_H$ and $(\id\otimes\pi)\circ\Delta_H$ respectively), while $C$ has a distinguished grouplike element $\ol 1$, where $1\in H$ is the unit. Write 
	\[
		{^\pi H}={^C H}\{h\in H\ |\ \ol{h_1}\otimes h_2=\ol 1\otimes h\},
	\]
	\[
		H^\pi=H^C=\{h\in H\ |\ h_1\otimes\ol{h_2}=h\otimes\ol 1\}.
	\]
These are what we were calling $\cat{LKer}(\pi)$ and $\cat{RKer}(\pi)$ back in the introduction, following the notation in \cite{Andruskiewitsch1995}. They are the spaces of $\ol 1$-coinvariants under the left and right coaction of $C$ on $H$ respectively, in the sense of \cite[28.4]{MR2012570}.

Dually, let $\iota:A\to H$ be an algebra map, and set $A^+=\iota^{-1}(\ker\varepsilon_H)$. Write $H_\iota=H_A$ for the left $H$-module $H/H\iota(A^+)$, and similarly, ${_\iota H}={_AH}=H/\iota(A^+)H$. 
	
It is now an easy exercise to check that if $\iota:A\to H$ is a right coideal subalgebra, then $H_A$ is a quotient left module coalgebra, and vice versa, if $\pi:H\to C$ is the projection on a quotient left module coalgebra, then $^C H$ is a right coideal subalgebra of $H$. 	
	\[
		\tikz[anchor=base]{
  		\path (0,0) node[text width=3.5cm] (algs) {set of right coideal subalgebras of $H$} +(6,0) node[text width=4.5cm] (coalgs) {set of quotient left module coalgebras of $H$};
  		\draw[->] (algs) .. controls (1,1) and (5,1) .. (coalgs) node[pos=.5,auto] {$A\mapsto H_A$};
  		\draw[<-] (algs) .. controls (1,-1) and (5,-1) .. (coalgs) node[pos=.5,auto,swap] {${^CH}\mapsfrom C$};
  	}
 	\]
are order-reversing maps with respect to the obvious poset structures on the two sets (whose partial orders we write as $\preceq$)

\begin{remark}
Note that the two order-reversing maps form a \define{Galois connection} in the sense of \cite[IV.5]{MR1712872} between the two posets of right coideal subalgebras and left module quotient coalgebras.
\end{remark}

\begin{definition}\label{def.reflections}
Let $\iota:A\to H$ be a right coideal subalgebra, and $\pi:H\to C$ a quotient left module coalgebra. We call $\pi:H\to H_A$ (or $H_A$ itself) the \define{right reflection} of $\iota:A\to H$ or of $A$, and $\iota:{^CH}\to H$ (or $^CH$ itself) the \define{left reflection} of $\pi:H\to C$. We also write $r(A)$ and $r(C)$ for $H_A$ and $^CH$.
\end{definition}

Using this language, recall from \cite[1.2.3]{Andruskiewitsch1995}:

\begin{definition}\label{def.exact}	
Let $H$ b a Hopf algebra. For a right coideal subalgebra $A\to H$ and a quotient left module coalgebra $H\to C$ we say that $k\to A\to H\to C\to k$ is \define{exact} if $A$ and $C$ they are each other's reflections. 
\end{definition}

We usually drop the $k$'s and talk just about exact sequences $A\to H\to C$. 

If $H$ is a Hopf algebra and $C$ is a left $H$-module coalgebra, then $_H^C\Mm$ will be the category of left $H$-modules endowed with a left $C$-comodule structure which is a left $H$-module map from $M$ to $C\otimes M$ (where the latter has the left $H$-module structure induced by the comultiplication on $H$). Similarly, if $A$ is a right $H$-comodule algebra, then $\Mm_A^H$ is the category of vector spaces right $H$-comodules with a right $A$-module structure such that $M\otimes A\to M$ is a map of right $H$-comodules. The mophisms in each of these categories are required to preserve both structures. 
	
	Let $\iota:A\to H$ be a right coideal subalgebra and $\pi:H\to C$ a quotient left module coalgebra such that $\pi\circ\iota$ factors through $A\ni a\mapsto\varepsilon(a)\ol 1\in C$ (this is equivalent to saying that $A\preceq r(C)$, or $C\preceq r(A)$, in the two posets discussed before \Cref{def.reflections}). Then, there is an adjunction between the categories $_A\Mm$ and $_H^C\Mm$, and dually, an adjunction between $\Mm_A^H$ and $\Mm^C$. We will recall briefly how these are defined, omitting most of the proofs, which are routine. 
	
	Let $M\in{_A\Mm}$. The vector space $H\otimes_AM$ then has a left $H$-module structure, as well as a left $C$-comodule structure inherited from the left $C$-coaction on $H$ (checking this is where the condition $A\preceq r(C)$ is needed). This defines a functor $L:{_A\Mm}\to{_H^C\Mm}$. To go in the other direction, for $N\in{_H^C\Mm}$, let
	\begin{equation}\label{eqn.defn of R}
		R(N)=\{n\in N\ |\ n_{\langle -1\rangle}\otimes n_{\langle 0\rangle}=\ol 1\otimes n\}. 
	\end{equation}
This defines a functor, and as the notation suggests, $L$ is a left adjoint to $R$. 

	For the other adjunction, given $M\in\Mm_A^H$, define $L'(M)=M/MA^+$. This is a functor (with the obvious definition on morphisms), and it is left adjoint to $R':\Mm^C\to\Mm_A^H$ defined by $R'(N)=N\square_CH$; the latter has a right $H$-comodule structure obtained by making $H$ coact on itself, as well as a right $A$-module structure obtained from the right $A$-action on $H$.
	
	Let us now focus on the adjunction $_A\Mm\tofrom{_H^C}\Mm$. In \cite{MR549940}, the same discussion is carried out in a slightly less general situation: the adjunction described above is considered in the case $A=r(C)$. On the other hand, we remark that when $C=r(A)$, the category $_H^C\Mm$ introduced above is nothing but the category of \define{descent data} for the ring extension $A\to H$. Recall (\cite[Proposition 25.4]{MR2012570}) that in our case, this would be the category $^{H\otimes_AH}\Mm$ of left comodules over the canonical $H$-coring $H\otimes_AH$ associated to the algebra extension $A\to H$. This means left $H$-modules $M$ with an appropriately coassociative and counital left $H$-module map $\rho:M\mapsto (H\otimes_AH)\otimes_HM\cong H\otimes_AM$.   

	The usual bijection
\[
	H\otimes H \cong H\otimes H \qquad h\otimes k \mapsto h_1\otimes h_2k 
\]
is easily seen to descend to a bijection $H\otimes_AH\cong r(A)\otimes H$. Hence, we see that a map $\rho$ as above is the same thing as a map $\psi:M\mapsto r(A)\otimes M$. The other properties of $\rho$, namely being a coassociative, counital, left $H$-module map, precisely translate to $\psi$ being coassociative, counital, and a left $H$-module map respectively. Taking into account this equivalence $_H^{r(A)}\Mm\simeq{^{H\otimes_AH}\Mm}$, the adjunction $(L,R):{_A\Mm}\tofrom{_H^{r(A)}\Mm}$ is an equivalence as soon as $H$ is right faithfully flat over $A$ (this is the faithfully flat descent theorem; see \cite[Theorem 3.8]{MR1466623}). 

Apart from faithful flatness, there are other known criteria which ensure $(L,R)$ is an equivalence. To state one of them, let us recall some notation from \cite{MR2210664}.

For a ring $A$, consider the contravariant endo-functor ${_A}C_A$ on the category of $A$-bimodules defined by ${_A}C_A(M)=\mathrm{Hom}(M,\mathbb{Q}/\mathbb{Z})$; these are homomorphisms of abelian groups, with the usual $A$-bimodule structure induced from that on $M$. Then, \cite[Theorem 8.1]{MR2210664} (very slightly rephrased) reads:

\begin{theoremN}
If $\iota:A\to H$ is a map of rings such that ${_A}C_A(\iota):{_A}C_A(H)\to{_A}C_A(A)$ is a split epimorphism, then $H\otimes_A$ is an equivalence between $_A\Mm$ and $^{H\otimes_AH}\Mm$. 
\end{theoremN}

Since we have just observed that in our case the functor $H\otimes_A$ from the statement of the theorem can be identified with $L:{_A}\Mm\to{_H^{r(A)}\Mm}$, we get the following result as a consequence:

\begin{proposition}\label{prop.split_comonad}
With the previous notations, $(L,R):{_A\Mm}\tofrom{_H^{r(A)}\Mm}$ is an equivalence if the inclusion $\iota:A\to H$ splits as an $A$-bimodule map.\qedhere  
\end{proposition}

\begin{remark}
The paper \cite{MR2210664} deals with rings rather than Hopf algebras. To deduce \Cref{prop.split_comonad} one first uses the noted identification $_H^{r(A)}\Mm\simeq{^{H\otimes_AH}\Mm}$ to turn the problem into the usual formulation of descent for arbitrary rings. \cite[$\S$ 7,8]{MR2210664} spell this out.  
\end{remark}

As a kind of converse to the faithfully flat descent theorem, $(L,R)$ being an equivalence implies that $H$ is right $A$-faithfully flat. Indeed, $H\otimes_A$ is then exact on $_A\Mm$. Note that we are using the fact that $_H^{r(A)}\Mm$ is abelian, with the same exact sequences as $\cat{Vec}$. All in all, this proves

\begin{proposition}\label{prop.equiv_iff_flat}
	Let $\iota: A\to H$ be a right coideal subalgebra. Then, the adjunction $(L,R):{_A\Mm}\tofrom{_H^{r(A)}\Mm}$ is an equivalence iff $H$ is right $A$-faithfully flat. \qedhere
\end{proposition}  
\begin{remark}
	This result is very similar in spirit to the equivalence (5) $\iff$ (3) in \cite[Theorem I]{MR1098988}, or to (1) $\iff$ (2) in \cite[Lemma 1.7]{MR2163407}. These can all be deduced from much more general, coring-flavored descent theorems that are now available, such as, say, \cite[Theorem 2.7]{MR2247888}.
\end{remark}

\subsection{CQG algebras}\label{subse.CQG}

For background, we rely mainly on \cite[11.3, 11.4]{MR1492989} or the paper \cite{MR1310296}, where these objects were originally introduced. Recall briefly that these Hopf algebras are meant to have just enough structure to imitate algebras of representative functions on compact groups. This means they are complex $*$-algebras (i.e. they possess conjugate-linear involutive multiplication-reversing automorphisms $*$) as well as Hopf algebras, and the two structures are compatible in the sense that the comultiplication and the counit are both $*$-algebra homomorphisms. 

In addition, CQG (Compact Quantum Group) algebras are required to have \define{unitarizable} comodules. This is a condition we will not spell out in any detail, but it says essentially that every finite-dimensional comodule has an inner product compatible with the coaction in some sense (once more imitating the familiar situation for compact groups, where invariant inner products on representations can be constructed by averaging against the Haar measure). In particular, CQG algebras are automatically cosemisimple, and hence fit comfortably into the setting of \Cref{se.main}.  

Not all $*$-algebras have enveloping $C^*$-algebras, but CQG algebras do. See, e.g. \cite[11.3.3]{MR1492989}. Such a completion is a so-called \define{full}, or \define{universal} $C^*$-algebraic compact quantum group, in the sense that it is a (unital) $C^*$-algebra $A$ endowed with coassociative $C^*$-algebra homomorphism $A\to A\otimes A$ (minimal $C^*$ tensor product) with additional conditions (\cite[11.3.3, Proposition 32]{MR1492989} or \cite[$\S$4,5]{MR1310296}). 

On the very few occasions when tensor product $C^*$-algebras come up, $\otimes$ always denotes the smallest $C^*$ tensor product (as treated in \cite[T.5]{MR1222415}, for instance). The term \define{completely positive map} between $C^*$-algebras will also make brief appearances. Recall that a linear map $T:A\to B$ between $C^*$-algebras is said to be positive if for each $x\in A$ we have $T(x^*x)=y^*y$ for some $y\in B$, and completely positive \cite[IV.3]{MR1873025} if the maps $\id\otimes T:M_n\otimes A\to M_n\otimes B$ between matrix algebras are all positive.


\section{Main results}\label{se.main}

We now prove the statement from the title of the paper:

\begin{theorem}\label{th.main}
A cosemisimple Hopf algebra is faithfully flat over all its Hopf subalgebras.
\end{theorem}
\begin{proof}
Let $H$ be cosemisimple, and $\iota:A\to H$ an inclusion of a Hopf subalgebra. Combining \Cref{prop.equiv_iff_flat,prop.split_comonad}, it suffices to show that $\iota$ splits as an $A$-bimodule map. In fact, one can even find a subcoalgebra $B\le H$ with $H=A\oplus B$ as $A$-bimodules.

Let $I$ be the set of simple subcoalgebras of $H$, and $J$ the subset of $I$ consisting of subcoalgebras contained in $A$. One then has $H=\bigoplus_ID$, and $A=\bigoplus_JD$. Define $B=\bigoplus_{I\setminus J}D$; in other words, $B$ is the direct sum of those simple subcoalgebras of $H$ which are not in $A$. Clearly, $B$ is a subcoalgebra, and $H=A\oplus B$, and we now only need to check that $B$ is invariant under (either left or right) multiplication by $A$. 
	
Let $D\in J$ and $E\in I\setminus J$ be simple subcoalgebras of $A$ and $B$ respectively. The product $ED$ inside $H$ is then $\cat{coalg}(V_E\otimes V_D)$ (see last paragraph above \cref{subse.descent}). Now assume $F\in J$ is a summand of $ED$. This means that $V_F\le V_E\otimes V_D$, so $V_E^*\le V_D\otimes V_F^*$. This is absurd: $V_E^*$ is a $B$-comodule, while $V_D\otimes V_F^*$ is an $A$-comodule. 
\end{proof}

\begin{remark}
This proves the first part of \cite[Conjecture 1]{MR2504527}; the second part, stating the faithful coflatness of a CQG algebra over quotient CQG algebras, follows immediately from the cosemsimplicity of CQG algebras. 
\end{remark}

\begin{remark}
Examples of cosemisimple Hopf algebras which are not faithfully coflat over \define{quotient} Hopf algebras abound, at least in characteristic zero. 

Indeed, let $G$ be a reductive complex algebraic group, and $B$ a Borel subgroup. Denoting by $\Oo(\bullet)$ `regular functions on the variety $\bullet$', the Hopf algebra $H=\Oo(G)$ is cosemisimple (e.g. \cite[p. 178]{MR0240104}, and it surjects onto $C=\Oo(B)$. 

If the surjection $H\to C$ were to be faithfully coflat, then, by \cite[Theorem 2]{MR549940}, we could reconstruct $C$ as $H/HA^+$ for $A=r(C)$. But $A$ is simply the algebra of global regular functions on the projective variety $G/B$, and hence consists only of constants; this provides the contradiction.    
\end{remark}

In fact, the result can be strengthened slightly. Recall that the \define{coradical} $C_0$ of a coalgebra $C$ is the sum of all its simple subcoalgebras.

\begin{corollary}\label{cor.main}
A Hopf algebra $H$ whose coradical $H_0$ is a Hopf subalgebra is faithfully flat over its cosemisimple Hopf subalgebras. 
\end{corollary}
\begin{proof}
Any cosemisimple Hopf subalgebra $A\le H$ will automatically be contained in the coradical $H_0$. By the previous corollary, $H_0$ is faithfully flat over $A$. On the other hand, Hopf algebras are faithfully flat (and indeed free) over sub-bialgebras which contain the coradical (\cite[Corollary 1]{MR0437582}); in particular, in this case, $H$ is faithfully flat over $H_0$. The conclusion follows. 
\end{proof}

Now let us place ourselves in the setting of \Cref{th.main}, assuming in addition that the Hopf subalgebra $A\to H$ is \define{conormal} in the language of \cite{Andruskiewitsch1995}. This simply means that $HA^+=A^+H$, and it is equivalent to $C=r(A)$ being a quotient Hopf algebra of $H$, rather than just a quotient coalgebra (\cite[Definition 1.1.9]{Andruskiewitsch1995}). Recalling the decomposition $H=A\oplus B$ as a direct sum of subcoalgebras, $C$ breaks up as the direct sum of the coalgebras $k=k\ol 1$ and $B/BA^+$. In other words, the coalgebra spanned by the unit of the Hopf algebra $C$ has a coalgebra complement in $C$. It follows (\cite[Theorem 14.0.3, (c) $\iff$ (f)]{MR0252485}) that $C$ is a cosemisimple Hopf algebra. Our aim, in the rest of this section, is to extend this result to the general case covered by \Cref{th.main}:

\begin{theorem}\label{th.mainbis}
If $\iota:A\to H$ is a Hopf subalgebra of a cosemisimple Hopf algebra $H$, then the coalgebra $C=r(A)$ is cosemisimple. 
\end{theorem}
\begin{proof}
We know from \Cref{th.main} that $H$ is right $A$-faithfully flat, and hence also left faithfully flat (just flip everything by means of the bijective antipode). This then implies, for example by \cite[Theorem 1]{MR549940}, that the second adjunction we introduced above, $(L',R'):{\Mm_A^H}\tofrom{\Mm^C}$ is an equivalence. It is then enough to show that all objects of the category $\Mm_A^H$ are projective, and this is precisely what the next two results do.	
\end{proof}

\begin{definition}
An object of $\Mm_A^H$ is said to be \define{$A$-projective} if it is projective as an $A$-module.
\end{definition}	

\begin{proposition}\label{prop.retract of Aproj}
Under the hypotheses of \cref{th.mainbis}, every object of $\Mm_A^H$ is $A$-projective. 
\end{proposition}
\begin{proof}
Let $M\in\Mm_A^H$ be an arbitrary object. Endow $M\otimes H$ with a right $H$-comodule structure by making $H$ coact on itself, and also a right $A$-module structure by the diagonal right action (i.e. $M\otimes H$ is the tensor product in the monoidal category $\Mm_A$). It is easy to check that these are compatible in the sense that they make $M\otimes H$ into an object of $\Mm_A^H$, and the map $\rho:m\mapsto m_{\langle 0\rangle}\otimes m_{\langle 1\rangle}\in M\otimes H$ giving $M$ its right $H$-comodule structure is actually a morphism in $\Mm_A^H$. Similarly, $\id\otimes\varepsilon_H:M\otimes H\to M$ is a morphism in $\Mm_A$, and it splits the inclusion $\rho$. It follows that it is enough to show that the object $M\otimes H\in\Mm_A^H$ described above is $A$-projective.
	
	\Cref{th.main} says that $H$ is $A$-faithfully flat, and it follows from \cite[Corollary 2.9]{MR1302835} that it is then (left and right) $A$-projective. This means that $M\otimes H$ can be split embedded (in the category $\Mm_A$) into a direct sum of copies of $M\otimes A$, with the diagonal right action of $A$. But
	\[
		M\otimes A \longto M\otimes A, \qquad  m\otimes a\mapsto ma_1\otimes a_2
	\]
exhibits an isomorphism from $M\otimes A$ with the right $A$-action on the right tensorand to $M\otimes A$ with the diagonal $A$-action (its inverse is $m\otimes a\mapsto mS(a_1)\otimes a_2$). This means that in $\Mm_A$, $M\otimes H$ is a direct summand of a direct sum of copies of $A$, i.e. projective. 	 
\end{proof}
		
\begin{proposition}\label{prop.Aproj means proj}	
	Under the hypotheses of \cref{th.mainbis}, $A$-projective objects of $\Mm_A^H$ are projective. 
\end{proposition}

	Before going into the proof, we need some preparation, including additional notation to keep track of the several $A$-module or $H$-comodule structures that might exist on the same object. 
	
	As in the proof of \Cref{th.main}, denote by $I$ and $J\subseteq J$ the sets of simple right comodules over $H$ and $A$, respectively. Recall that these are also in one-to-one correspondence with the simple subcoalgebras of $H$ and $A$, respectively. We will henceforth denote by $\varphi:H\to A$ the map which is the identity on $A$, and sends every simple subcoalgebra $D\in I\setminus J$ to $0$.
		
	Notice now that $A$ acts on $H$ (as well as on itself) not just by the usual right regular action, but also by the right adjoint action: $h\triangleleft a=S(a_1)ha_2$ ($h\in H$, $a\in A$). This gives $H$ and $A$ a second structure as objects in $\Mm_A^H$. When working with this structure rather than the obvious one, we denote these objects by $H_{\ad}$ and $A_{\ad}$. 
	
\begin{lemma}\label{lemma.varphi}
	(a) For any object $M\in\Mm_A^H$, $M\otimes H_{\ad}$ becomes an object of $\Mm_A^H$ when endowed with the diagonal $A$-action (where $A$ acts on $M\in\Mm_A^H$ and on $H$ by the right adjoint action) and the diagonal $H$-coaction.   
	
	(b) Similarly, $M\otimes A_{\ad}\in\Mm_A^H$.
	
	(c) $\id\otimes\varphi:M\otimes H_{\ad}\to M\otimes A_{\ad}$ respects the structures from (a) and (b), and hence is a morphism in $\Mm_A^H$. 
\end{lemma}
\begin{proof}
	We will only prove (a); (b) is entirely analogous, while (c) follows immediately, since $\varphi$ clearly preserves both the right $H$-coaction and the adjoint $A$-action. 
	
	Proving (a) amounts to checking that the diagram
	\[
		\tikz[anchor=base]{
  		\path (0,0) node (1) {$M\otimes H_{\ad}\otimes A$} +(4,0) node (2) {$M\otimes H_{\ad}$} +(4,-2) node (3) {$M\otimes H_{\ad}\otimes H$} +(0,-2) node (4) {$M\otimes H_{\ad}\otimes H\otimes A$};
  		\draw[->] (1) -- (2);
  		\draw[->] (4) -- (3);
  		\draw[->] (1) -- (4);
  		\draw[->] (2) -- (3); 
 		}
 	\]
is commutative. The path passing through the upper horizontal line is
	\[
		m\otimes h\otimes a\longmapsto ma_1\otimes S(a_2)ha_3\longmapsto m_0a_1\otimes S(a_4)h_1a_5\otimes m_1a_2S(a_3)h_2a_6,
	\]
while the other composition is 
	\[
		m\otimes h\otimes a\longmapsto m_0\otimes h_1\otimes m_1h_2\otimes a\longmapsto m_0a_1\otimes S(a_2)h_1a_3\otimes m_1h_2a_4.  
	\] 	
Using the properties of the antipode and counit in a Hopf algebra, we have
	\begin{alignat*}{-1}
		m_0a_1\otimes S(a_4)h_1a_5\otimes m_1a_2S(a_3)h_2a_6 &=& &~m_0a_1\otimes S(\varepsilon(a_2)a_3)h_1a_4\otimes m_1h_2a_5\\ &=& &~m_0a_1\otimes S(a_2)h_1a_3\otimes m_1h_2a_4,
	\end{alignat*}
concluding the proof	
\end{proof}

	Now denote by $(M\otimes H)^r\in\Mm_A^H$ the object from the proof of \cref{prop.retract of Aproj}: the $A$-action is diagonal, while $H$ coacts on the right tensorand alone. The upper $r$ is meant to remind the reader of this. 
	
\begin{lemma}\label{lemma.psi}
	For $M\in \Mm_A^H$, $\psi_M:M\otimes H\to M\otimes H$ defined by
	\[
		m\otimes h\longmapsto m_0\otimes S(m_1)h. 
	\]
is a morphism in $\Mm_A^H$ from $(M\otimes H)^r$ to $M\otimes H_{\ad}$. 
\end{lemma} 
\begin{proof}
	We only check compatibility with the $A$-actions, leaving the task of doing the same for $H$-coactions to the reader. The composition $(M\otimes H)^r\otimes A\longto (M\otimes H)^r\longto M\otimes H_{\ad}$ is
	\[
		m\otimes h\otimes a\verylongmapsto{} ma_1\otimes ha_2\verylongmapsto{\psi_M} m_0a_1\otimes S(m_1a_2)ha_3, 
	\]
while the other relevant composition is
	\[
		m\otimes h\otimes a\verylongmapsto{\psi_M\otimes\id} m_0\otimes S(m_1)h\otimes a\verylongmapsto{} m_0a_1\otimes S(a_2)S(m_1)ha_3.
	\]
Since $S$ is an algebra anti-morphism, they are equal. 
\end{proof}

	Finally, we have
	
\begin{lemma}\label{lemma.mult}
	Let $M\in\Mm_A^H$. The map $M\otimes A\to M$ giving $M$ its $A$-module structure is a morphism $M\otimes A_{\ad}\to M$ in $\Mm_A^H$.  
\end{lemma} 
\begin{proof}
	Compatibility with the $H$-coactions is built into the definition of the category $\Mm_A^H$, so one only needs to check that the map is a morphism of $A$-modules. In other words, we must show that the diagram 
	\[
		\tikz[anchor=base]{
  		\path (0,0) node (1) {$M\otimes A_{\ad}\otimes A$} +(4,0) node (2) {$M\otimes A_{\ad}$} +(4,-2) node (3) {$M$} +(0,-2) node (4) {$M\otimes A$};
  		\draw[->] (1) -- (2);
  		\draw[->] (4) -- (3);
  		\draw[->] (1) -- (4);
  		\draw[->] (2) -- (3); 
 		}
 	\]
is commutative. The right-down composition is
	\[
		m\otimes a\otimes b\verylongmapsto{} mb_1\otimes S(b_2)ab_3\verylongmapsto{} mb_1S(b_2)ab_3,
	\]
while the other composition is 
	\[
		m\otimes a\otimes b\verylongmapsto{} ma\otimes b \verylongmapsto{} mab;
	\]
they are thus equal. 
\end{proof}

\begin{lemma}\label{lemma.split}
	For $M\in\Mm_A^H$, the composition
	\[
		\tikz[anchor=base]{
  		\path (0,0) node (1) {$t_M:(M\otimes H)^r$} +(4,0) node (2) {$M\otimes H_{\ad}$} +(8,0) node (3) {$M\otimes A_{\ad}$} +(10,0) node (4) {$M$};
  		\draw[->] (1) -- (2) node[pos=.5,auto] {$\scriptstyle \psi_M$};
  		\draw[->] (2) -- (3) node[pos=.5,auto] {$\scriptstyle \id\otimes\varphi$};
  		\draw[->] (3) -- (4); 
 		}
	\]
where the last arrow gives $M$ its $A$-module structure is a natural transformation from the $\Mm_A^H$-endofunctor $(\bullet\otimes H)^r$ to the identity functor, and it exhibits the latter as a direct summand of the former.
\end{lemma}
\begin{proof}
	The fact that $t_M$ is a map in $\Mm_A^H$ follows from \cref{lemma.varphi,lemma.psi,lemma.mult}. Naturality is immediate (one simply checks that it holds for each of the three maps), as is the fact that $t_M$ is a left inverse of the map $M\to(M\otimes H)^r$ giving $M$ its $H$-comodule structure. 
\end{proof}

	We are now ready to prove the result we were after.
	
\begin{proof of Aproj}
	Let $P\in\Mm_A^H$ be an $A$-projective object. We must show that $\Mm_A^H(P,\bullet)$ is an exact functor. Embedding the identity functor as a direct summand into $(\bullet\otimes H)^r$ (\cref{lemma.split}), it suffices to show that $\Mm_A^H(P,(\bullet\otimes H)^r)$ is exact. 
	
	$(\bullet\otimes H)^r:\Mm_A\to\Mm_A^H$ is right adjoint to $\text{forget}:\Mm_A^H\to\Mm_A$ (as $\Mm_A^H$ is the category of coalgebras for the comonad $\bullet\otimes H$ on $\Mm_A$; \cite[Theorem VI.2.1]{MR1712872}), so $\Mm_A^H(P,(\bullet\otimes H)^r)$ is naturally isomorphic to $\Mm_A(P,\bullet)$, which is exact by our assumption that $P$ is $A$-projective.
\end{proof of Aproj}

\begin{remark}
	In the above proof, the forgetful functor $\text{forget}:\Mm_A^H\to\Mm_A$ has been suppressed in several places, in order to streamline the notation; we trust that this has not caused any confusion.  
\end{remark}

\begin{remark}
	The proof of \cref{prop.retract of Aproj} is essentially a rephrasing of the usual proof that Hopf algebras $H$ with a (right, say) integral sending $1_H$ to $1$ are cosemisimple (\cite[$\S$14.0]{MR0252485}; we will call such integrals unital). The map $\varphi:H\to A$ introduced in \cref{lemma.varphi} might be referred to as an \define{$A$-valued right integral} (by which we mean a map preserving both the right $H$-comodule structure and the right adjoint action of $A$), and specializes to a unital integral when $A=k$. In conclusion, one way of stating \cref{prop.Aproj means proj} would be:
	
	\define{If the inclusion $\iota:A\to H$ of a right coideal subalgebra is split by an $A$-valued right integral, then the forgetful functor $\Mm_A^H\to\Mm_A$ reflects projectives.}	 
\end{remark}

\begin{remark}	
	\Cref{prop.retract of Aproj,prop.Aproj means proj} can both be traced back to work by Y. Doi, but we have included proofs for completeness. \Cref{prop.retract of Aproj}, for instance, is a consequence of \cite[Theorem 4]{MR688207}. Similarly, \Cref{prop.Aproj means proj} follows from \cite[Theorem 1]{MR1098982}. I thank the referee for pointing this out. 
\end{remark}


\section{Expectations on CQG subalgebras are positive}\label{se.CQG}

We now move the entire $A\to H\to C$ setting over to the case when $H$ is a CQG algebra. We take for granted the preceding sections, and in particular the fact that $C$ is cosemisimple (\Cref{th.mainbis}). The inclusion $\iota: A\to H$ is now one of $*$-algebras, and we follow operator algebraists' convention of referring to its left inverse $p:H\to A$ from the proof of \Cref{th.main} as the \define{expectation} of $H$ on $A$ (in accordance with a view of $A$ and $H$ as consisting of random variables on non-commutative measure spaces). Positivity here means the following:

Think of $H$ as embedded in its universal $C^*$ completion $H_u$ (\Cref{subse.CQG}), and complete $A$ to $A_u$ with the subspace norm. Then, $p$ extends to a completely positive map $H_u\to A_u$. Equivalently, the self-map $\iota\circ p:H\to H$ lifts to a completely positive self-map of $H_u$.

Note that a functional $\psi\in H^*$ with $\psi(1)=1$ extends to a state on the $C^*$ completion $H_u$ if and only if it is positive in the usual sense, i.e. $\psi(x^*x)\ge 0$ for every $x\in H$.  

The main result of the section is:

\begin{theorem}\label{th.CQG}
Let $\iota:A\to H$ be an inclusion of CQG algebras. Then, the expectation $p:H\to A$ is positive in the above sense. 
\end{theorem}

\begin{remark}
So called \define{expected} $C^*$-subalgebras of (locally) compact quantum groups have featured prominently in the literature (see \cite{MR2335776,MR2999995} and references therein). The techniques used in the proof of \Cref{th.CQG} will be applied elsewhere \cite{C} to characterize all right coideal $*$-subalgebras $A$ of a CQG algebra $H$ which are expected in the sense of admitting a positive splitting of the inclusion as an $A$-bimodule, right $H$-comodule map, where positivity is understood as in \Cref{th.CQG}.  
\end{remark}

Let us first reformulate the theorem slightly. Denote the unique unital (left and right) integral of $C$ by $h_C$, and the composition $h_C\circ\pi$ by $\varphi$ (where $\pi:H\to C$ is the surjection we start out with). The expectation decomposes as $(\varphi\otimes\id)\circ\Delta:H\to A$. This follows easily from the decomposition $H=A\oplus B$ as a direct sum of subcoalgebras used in the proof of \Cref{th.main} and the fact that $\varphi|_A$ equals $\varepsilon_A$ and $\varphi|_B$ is the zero map.

\begin{remark}\label{re.phi_preserves_*}
Let us note in passing that $\varphi$ is self-adjoint as a functional, in the sense that $\varphi(x^*)$ is the complex conjugate of $\varphi(x)$ for any $x\in H$. This follows immediately from $\varphi|_A=\varepsilon_A$ and $\varphi|_B=0$, the fact that $A$ and $B$ are closed under $*$, and the fact that $\varepsilon$ is a $*$-homomorphism.

This observation is needed in the proof of item 2. in \Cref{prop.F}, for instance. 
\end{remark}

\begin{lemma}\label{lemma.idempotent}
The conclusion of \Cref{th.CQG} holds if and only if the functional $\varphi\in H^*$ is positive. 
\end{lemma}
\begin{proof}
Note that $\varphi$ equals $\varepsilon\circ p$ (more pedantically, in this expression $\varepsilon$ is the restriction of $\varepsilon_H$ to $A$). If $p$ is positive then so will $\varphi$, given that $\varepsilon$ is a $*$-algebra map $A\to \CC$ which lifts to $A_u$. 

Conversely, if $\varphi$ is positive (and hence lifts to a state on $H_u$), then both maps in the composition $(\varphi\otimes\id)\circ\Delta:H\to H$ lift to completely positive maps on the apropriate $C^*$ completions ($\Delta$ lifts to a $C^*$-algebra map $H_u\to H_u\otimes H_u$, while $\varphi\otimes\id:H_u\otimes H_u\to H_u$ will also be completely positive). But that composition is precisely $\iota\circ p$, as noted above. 
\end{proof}

\begin{remark}\label{rem.phi_S-inv}
The identity $\varphi=\varepsilon\circ p$, in particular, shows that $\varphi\circ S=\varphi$. This is needed below.  
\end{remark}

We are going to take what looks like a detour to make the necessary preparations. 

For a cosemisimple coalgebra $D$ over an algebraically closed field, denote by $D^\bullet$ its \define{restricted dual}: It is the direct sum of the matrix algebras dual to the matrix subcoalgebras of $D$. In general, $D^\bullet$ is a non-unital algebra. In our case, the full dual $H^*$ is in addition a (unital) $*$-algebra, with $*$ operation defined by
\begin{equation}\label{eq.dual_star}
	f^*(x) = f((Sx)^*)^*,\ \forall x\in H,
\end{equation}  
where the outer $*$ means complex conjugation of a number (see e.g. \cite[4.3]{MR1658585}). Furthermore, $C^\bullet\le H^*$ is a $*$-subalgebra. 

Finally, again for a cosemisimple coalgebra $D$, we will talk about its completion $\ol{D}$; it is by definition the direct product of the matrix subcoalgebras comprising $D$. Equivalently, $\ol{D}$ is the (ordinary, vector space) dual of $D^\bullet$. The module structure $H\otimes C\to C$ extends to an action of $H$ on $\ol{C}$.

\begin{remark}\label{rem.mod_cat}
This extension of the $H$-module structure to $\ol{C}$ is a simple enough observation, but there is some content to it. The claim is that for $x\in H$ and some simple subcoalgebra $C_\alpha\le C$ (for $\alpha\in\widehat{C}$), there are only finitely many simples $\beta\in\widehat{C}$ such that $xC_\beta$ intersects $C_\alpha$ non-trivially. 

Although $\Mm^C$ is not monoidal, $V\otimes W$ can be made sense of as a $C$-comodule for any $H$-comodule $V$ and $C$-comodule $W$. This makes $\Mm^C$ into a \define{module category} over the monoidal category $\Mm^H$. Upon rephrasing the claim using the correspondence $W\mapsto \cat{coalg}(W)$ between comodules and subcoalgebras, it reads: For each finite-dimensional $H$-comodule $V$ and each $\alpha\in\widehat{C}$, there are only finitely many $\beta\in\widehat{C}$ such that (identifying $\alpha,\beta$ with the corresponding comodules) $\Hom_{\Mm^C}(\alpha,V\otimes\beta)$ is non-zero. But just as in a rigid monoidal category, $V\otimes:\Mm^C\to\Mm^C$ is right adjoint to $V^*\otimes$, and hence we're saying only finitely many $\beta$ satisfy $\Hom(V^*\otimes\alpha,\beta)\ne 0$. This is clear simply because $V^*\otimes\alpha$ is some finite direct sum of irreducibles.   
\end{remark}

First, a preliminary result:

\begin{lemma}\label{lemma.S2pos}
The squared antipode $S^2$ of $H$ descends to an automorphism of every simple subcoalgebra $C_\alpha$ of $C$. Moreover, the resulting automorphism on the $C^*$-algebra $C_\alpha^*$ is conjugation by an invertible positive operator. 
\end{lemma}
\begin{proof}
That $S^2$ descends to $C=H/HA^+$ is clear from the fact that it acts on $A$. We move the action over to duals by precomposition: $S^2f = f(S^2\cdot)$ for $f\in H^*$. 

Now let $D\le H$ be a simple subcoalgebra, and $\bigoplus_{\alpha\in I} C_\alpha$, $I\subset\widehat{C}$ the image of $D$ through $H\to C$. The squared antipode acts on $D^*$ as conjugation by a positive operator $F$ (\cite[11.30, 11.34]{MR1492989}), and, by the previous paragraph, preserving the subalgebra $B=\bigoplus_{\alpha\in I} C_\alpha^*$. In particular, conjugation by $F$ permutes the $|I|$ minimal non-zero projections $p_\alpha$, $\alpha\in I$ in the center of $B$. I claim that this permutation action is in fact trivial, which would finish the proof. 

To check the claim, consider the unique (up to isomorphism) simple $*$-representation of $D^*$ on a Hilbert space $\mathcal H$. If $Fp_\alpha F^{-1}$ were equal to some $p_\beta$ with $\beta\ne\alpha\in I$, then $F$ would map the range of $p_\alpha$ onto the range of $p_\beta$. Denoting by $\langle\ ,\ \rangle$ the inner product on $\mathcal{H}$, this implies that $\langle Fx,x\rangle$ vanishes for any $x$ in the range of $p_\alpha$. This cannot happen for non-zero $x$, as $F$ is both positive and invertible. 
\end{proof}

We now establish the existence of a kind of ``relative Haar measure'' on $C^\bullet$.

\begin{proposition}\label{prop.rel_haar}
There is an element $\theta\in \ol C$ satisfying the following conditions:
\renewcommand{\labelenumi}{(\alph{enumi})}
\begin{enumerate}
	\item Writing $\theta$ as a formal sum of elements in the simple subcoalgebras of $C$, its component in $\CC\ol 1\le C$ is $\ol 1$; 
	\item It is $H$-invariant, in the sense that $x\theta=\varepsilon(x)\theta$ for $x\in H$;
	\item It is positive as a functional on the $*$-algebra $C^\bullet$.
\end{enumerate}
\end{proposition}
\begin{sketch}
Let $e_i$, $e^i$, $i\in I$ be dual bases in $C$ and $C^\bullet$ respectively, compatible with the decomposition of $C$ into simple subcoalgebras. We distinguish an element $0\in I$ such that $e_0=\ol 1$. Since the automorphism $S^2$ of $H$ descends to $C=H/HA^+$, the definition
\[
	\theta = \sum_{i\in I} e^i(S^2 e_{i(2)}) e_{i(1)}
\]
makes sense as an element of $\ol C$, and clearly satisfies (a). Moreover, the definition does not depend on the choice of bases. 

The calculation proving $H$-invariance can simply be lifted e.g. from \cite[1.1]{MR1415374}. Even though that result is about finite-dimensional Hopf algebras, it works verbatim in the present setting. 

Finally, let us prove positivity, this time imitating \cite{MR1378538}. Let $\alpha\in\widehat{C}$, and $u\in C_\alpha^*\le C^\bullet$ an element. We can assume harmlessly that the bases $e_i$, $e^i$ are organized as matrix (co)units, i.e. those $e_i$ in the matrix coalgebra $C_\alpha$ form a matrix counit $e_{pq}$, and $e^i$ will then be the dual matrix unit $e^{pq}\in C_\alpha^*$. 

Now note that $e_{pq}$, regarded as a functional on $C_\alpha^*$, can be written as $\tr_\alpha(\cdot\ e^{qp})$, where $\tr_\alpha$ is the trace on the matrix algebra $C_\alpha^*\cong M_n$, so that $\tr_\alpha(1)=n$. In conclusion, the component of $\theta$ in $C_\alpha$, regarded as a functional on $C_\alpha^*$, is 
\begin{equation}\label{eq.theta}
	\theta_\alpha = \sum_{p,q} \tr_\alpha(\cdot\ S^2(e^{pq})e^{qp}).
\end{equation}
If $Q\in C_\alpha^*$ is a positive operator such that conjugation by $Q$ equals $S^2$ on $C_\alpha^*$ (\Cref{lemma.S2pos}), then, suppressing summation over $p,q=\ol{1,n}$,
\[
	S^2(e^{pq})e^{qp} = Qe^{pq}Q^{-1}e^{qp} = \tr_\alpha(Q^{-1})Q.
\]
This is a positive operator, and the conclusion follows.
\end{sketch}

\begin{remark}\label{rem.theta_S2-inv}
The expression \Cref{eq.theta}, the invariance of $\theta$ with respect to bases, and the fact that $S^2(e^{pq})$ are again matrix units make it clear that $\theta\circ S^2=\theta$. In fact $\theta$ is unique, but we do not need this stronger fact.  
\end{remark}

\begin{definition}\label{def.relF}
Keeping the previous notations, the \define{$\varphi$-relative Fourier transform} $\Ff:H\to C^\bullet$ is defined as
\[
	H\ni x\mapsto \varphi(Sx\cdot). 
\]
\end{definition}

There is a slight abuse of notation in the definition: although a priori $\varphi$ is a functional on $H$, it descends to one on $C=H/HA^+$. The map $\Ff$ is a relative analogue to the usual Fourier transform \cite[$\S$2]{MR1059324}, and enjoys similar properties. Let us record some of them:

\begin{proposition}\label{prop.F}
The map $\Ff:H\to C^\bullet$ introduced above satisfies the following relations:
\begin{enumerate}
	\item $\Ff(x\triangleleft \Ff y) = \Ff x\Ff y$, $\forall x,y\in H$, where the right action $\triangleleft$ of $H^*$ on $H$ is defined by $x\triangleleft f = f(x_1)x_2$. 
	\item $\Ff(x)^* = S^2\Ff((Sx)^*)$, where the $*$ structure on $C^\bullet$ is defined in \Cref{eq.dual_star}, and $S^2$ acts on $H^*$ by precomposition, as in the proof of \Cref{lemma.S2pos}.
	\item $\varepsilon(x\triangleleft\Ff y) = \varphi(Sy\ x)$. 
	\item $\theta\Ff = \varepsilon$, where $\theta$ is the functional on $C^\bullet$ from \Cref{prop.rel_haar}. 
	\item $\Ff S^2=S^{-2}\Ff$.  
\end{enumerate}
\end{proposition} 
\begin{proof}
Most of this consists of simple computations, so let us only prove the first and fourth items. 

Applying both sides of 1. to $z\in H$, we have to prove \[ \varphi(Sy\ x_1)\varphi(Sx_2\ z) = \varphi(Sx\ z_1)\varphi(Sy\ z_2). \] Substituting $y$ for $Sy$, $z$ for $Sz$, and using $\varphi\circ S=\varphi$ (\Cref{rem.phi_S-inv}), this turns into \[ \varphi(yx_1)\varphi(zx_2) = \varphi(z_2x)\varphi(ySz_1). \] 

Now make the substitution $yx_1\otimes x_2=a\otimes b$, which in turn means $y\otimes x=aSb_1\otimes b_2$. The target identity turns into \[ \varphi(a)\varphi(zb)=\varphi(aSb_1Sz_1)\varphi(z_2b_2). \] Writing $zb=c$, it transforms further into \[ \varphi(a)\varphi(c)=\varphi(aSc_1)\varphi(c_2). \] Finally, the substitution of $c$ for $S^{-1}c$ and again $\varphi\circ S=\varphi$ turn this into \[ \varphi(a)\varphi(c) = \varphi(ac_2)\varphi(c_1). \] To prove this last equality, it suffices to split into two cases, according to whether $c$ is in $A$ or the complementary $A$-bimodule, right $H$-comodule $\ker(p)$. 

In the latter case, both $\varphi(c)$ and $\varphi(c_1)$ vanish. In the former, the left hand side is $\varphi(a)\varepsilon(c)$, while the right hand side is $\varphi(ac)$ (since $\varphi(c_1)=\varepsilon(c_1)$). These two expressions are equal because $\varphi=\varepsilon p$, and $p$ is an $A$-bimodule map.  

We now check 4. Applying its left hand side to $x\in H$, we get $\theta(\varphi(Sx\cdot))=\varphi(Sx\ \theta)$, where this time $\theta$ is thought of as an element of $\ol{C}$, $Sx\ \theta$ is the action of $Sx$ on it (\Cref{rem.mod_cat}), and $\varphi$ is regarded naturally as a functional on $\ol{C}$. By the $H$-invariance of $\theta$ (\Cref{prop.rel_haar} (b)), the expression is $\varepsilon(x)\varphi(\theta)=\varepsilon(x)$ by \Cref{prop.rel_haar} (a).  
\end{proof}

All of the ingredients are now in place.

\begin{proof of thCQG}
According to \Cref{lemma.idempotent}, it suffices to show that $\varphi(x^*x)\ge 0$ for all $x\in H$. We do this through a string of equalities based on the preliminary results of this section. 

Let $x,y\in H$. Then, we have
\begin{multline*}
	\theta((\Ff y)^*\Ff x)\ \overset{2.}{=}\ \theta(S^2\Ff((Sy)^*)\Ff x)\ =\ \theta(S^2\Ff((Sy)^*)\Ff (S^2 x))\ \overset{1.}{=}\ \theta\Ff((Sy)^*\triangleleft\Ff(S^2 x)) \\ 
	\overset{4.}{=}\ \varepsilon((Sy)^*\triangleleft\Ff(S^2 x))\ \overset{3.}{=}\ \varphi(S^3 x(Sy)^*)\ =\ \varphi(S^3 x S(S^2y)^*)\ =\ \varphi((S^2 y)^*S^2x),  
\end{multline*}
where the numbers above the equal signs refer to the items in \Cref{prop.F}, the second equality follows from 5. and the fact that $\theta S^2=\theta$ (\Cref{rem.theta_S2-inv}), the next-to-last one is a simple manipulation valid in any Hopf $*$-algebra, and the last equality is based on $\varphi S=\varphi$ (\Cref{rem.phi_S-inv}). Since the left hand side is non-negative when $x=y$, so is the right hand side. This concludes the proof of the theorem. 
\end{proof of thCQG}

\begin{remark}\label{rem.Planch}
The equality obtained in the course of the proof should be thought of as a Plancherel theorem \cite[7.9]{MR1157815}, to the effect that the relative Fourier transform is an isometry with respect to the ``inner products'' induced by $\varphi$ and $\theta$. 
\end{remark}


\bibliographystyle{plain}

\end{document}